\documentclass[12pt,english,twoside,a4paper,reqno]{amsart}
\usepackage[T1]{fontenc}
\usepackage[latin9]{inputenc}
\usepackage{geometry}
\usepackage{float}
\usepackage{amsbsy}
\usepackage{amstext}
\usepackage{amsthm}
\usepackage{amssymb}
\usepackage{graphicx}
\usepackage{color}
\usepackage{esint}
\usepackage{enumerate}
\PassOptionsToPackage{normalem}{ulem}
\usepackage{ulem}

\usepackage{babel}
\begin{document}

\title{Quantum Graphs via Exercises}

\author[R.~Band \& S.~Gnutzmann]{Ram Band$^{1}$and Sven Gnutzmann$^{2}$}

\address{$^{1}${\small{}Department of Mathematics, Technion--Israel Institute
of Technology, Haifa 32000, Israel}}

\address{$^{2}$School of Mathematical Sciences,   University of Nottingham,
Nottingham NG7 2RD, UK}
\begin{abstract}
Studying the spectral theory of Schr\"odinger operator on metric
graphs (also known as \textquotedbl{}quantum graphs\textquotedbl{})
is advantageous on its own as well as to demonstrate key concepts
of general spectral theory. There are some excellent references for
this study such as \cite{BerKuc_graphs} (mathematically oriented
book), \cite{GnuSmi_ap06} (review with applications to theoretical
physics), and \cite{Berkolaiko_arXiv17} (elementary lecture notes).
Here, we provide a set of questions and exercises which can accompany
the reading of these references or an elementary course on quantum
graphs. The exercises are taken from courses on quantum graphs which
were taught by the authors.
\end{abstract}

\maketitle

\section{Basic Spectral Theory of Quantum Graphs}

\subsection{$\delta$-type vertex conditions}
\label{subsec:Delta_type_conditions}

~\\[0.4cm]
\noindent
\textbf{Background}:
The most common vertex conditions of a graph are called Dirichlet
and Neumann. Note that the Neumann conditions are also sometimes called
Kirchhoff, standard or natural. Nevertheless, we stick to the name
Neumann in this note.

For the following set of questions we
consider a particle in a box described by the Schr\"odinger equation
\begin{equation}
-\frac{d^{2}f(x)}{dx^{2}}+V(x)f(x)=k^{2}f(x) \quad \textrm{ for } \, x\in(-l_1,l_2), \label{eq:schroedinger}
\end{equation}
where $l_1>0,\,l_2>0$ and $V(x)$ is called the electric potential and can be taken for example in $L^2([-l_1,l_2])$.
We employ here the Dirichlet boundary conditions $f(-l_{1})=0=f(l_{2})$, which physically mean that the particle is trapped in an infinite well, whose walls are at $x=-l_1$ and $x=l_2$. \\[0.2cm]

\noindent
\textbf{Questions:}\\[0.2cm]
\begin{itemize}
\item[\textbf{1.}]
  In this question, we choose the potential in (\ref{eq:schroedinger}) to be $V(x)=\alpha\delta(x)$, where $\alpha\in\mathbb{R}$ and $\delta(x)$ is Dirac's delta function. By doing so, we extend the original validity of equation (\ref{eq:schroedinger}), since the potential is now considered as a distribution and not merely a function. In order for the action of the $\delta$ distribution to be well-defined, we must require that $f$ is continuous at
  $x=0$ with $f_{0}:=f(0^{+})=f(0^{-})$.

  Show that any eigenfunction $f$ in (\ref{eq:schroedinger}) satisfies the
  matching condition
  \begin{equation}
  \alpha f_{0}=f'(0^{+})-f'(0^{-}). \label{eq:basic-matching-condition}
  \end{equation}
\\

  \textit{Comment}: Note that the sign of $\alpha$ determines whether the $\delta$-potential
is attractive ($\alpha<0$), repulsive ($\alpha>0$) or vanishing
($\alpha=0$). \\
 The relation (\ref{eq:basic-matching-condition}) is called the $\delta$-type
  vertex condition. Indeed, in our example $x=0$ may be considered as a graph vertex whose degree is two. If the vertex under consideration is of degree higher
  than two, then this condition is generalized and takes the form
  \begin{equation}
    af_{0}=\sum_{e}\left.f\right|_{e}'(0),\label{eq:Delta_type_conditions_general}
  \end{equation}
where the sum is taken over all edges adjacent to the considered vertex, and in all the summands the derivative is taken to be directed towards the vertex.
  In addition, we always require the continuity of the function at the
  considered vertex, i.e., all adjacent edges agree in the value they
  obtain at the vertex, and this value equals $f_0$. \\
  Note that the Dirichlet and Neumann conditions are obtained as special cases of the $\delta$-type vertex condition for some particular values of the parameter ($\alpha=0$ for Neumann and $\alpha\rightarrow\infty$ for Dirichlet).\\[0.2cm]
%\end{itemize}
%
%\pagebreak{}
%
%\noindent
%\textbf{Question:}
%\begin{itemize}
\item[\textbf{2.}]
  Next, we find the spectrum of the problem above:
  the interval $\left[-l_{1},l_{2}\right]$ with Dirichlet conditions
  at the end points,
  \begin{equation}
    f(-l_{1})=f(l_{2})=0,\label{eq:Dirichlet_conditions_of_interval}
  \end{equation}
  and a $\delta$-type vertex condition at $x=0$,
  \begin{align}
    f(0^{+}) & =f(0^{-})\label{Delta_type_conditions_1}\\
    \alpha f_{0}= & f'(0^{+})-f'(0^{-}).\label{eq:Delta_type_conditions_2}
  \end{align}
  ~\\[0.1cm]
  \begin{enumerate}[(a)]
  \item Show that the non-negative eigenvalues $k^{2}$ of the one-dimensional
    Laplacian, $-\frac{\rm{d}^2}{\rm{d}x^2}$, with the vertex conditions (\ref{eq:Dirichlet_conditions_of_interval}),(\ref{Delta_type_conditions_1}),(\ref{eq:Delta_type_conditions_2}) are given as the zeros of the following
    function
    \begin{equation}
      \zeta_{\alpha}(k)=\alpha+k\cot(kl_{1})+k\cot(kl_{2}).\label{eq:Secular_function_with_delta}
    \end{equation}
    We call such a function, whose zeros provide the graph's eigenvalues,
    a secular function.\\[0.1cm]
  \item For an irrational lengths ratio ($l_{1}/l_{2}\notin\mathbb{Q}$) show
    that there cannot be any eigenfunctions with $f_{0}=0$.

    Conclude that the poles of the secular function, (\ref{eq:Secular_function_with_delta}),
    cannot belong to the spectrum.\\[0.1cm]
  \item Show that $\zeta_{\alpha}(k)$ has only single poles if $l_{1}/l_{2}\notin\mathbb{Q}$.

    However, if $l_{1}/l_{2}\in\mathbb{Q}$ there are single and double
    poles of the secular function. Locate them and show that the double
    poles belong to the spectrum, while the single poles do not belong
    to it. This justifies the need of regularization of the secular function
    \begin{equation}
    \tilde{\zeta}_{\alpha}(k)=\zeta_{\alpha}(k)\sin(kl_{1})\sin(kl_{2}).\label{eq:secular-function-with-delta-no-poles}
    \end{equation}
    Namely, the zeros of the regularized secular function, $\tilde{\zeta}_{\alpha}(k)$,
    correspond to the Laplacian eigenvalues (irrespective to the value
    of $l_{1}/l_{2}$).\\[0.1cm]
  \item Consider an attractive $\delta$-potential ($\alpha<0$). Are there
    negative eigenvalues among the eigenvalues of the Schrödinger operator?
    Under what conditions and how many?\\[0.1cm]
  \item Show that if $\alpha=0$ in the problem above, this is equivalent
    to solving the eigenvalue problem of the Laplacian on the interval,
    $\left[-l_{1},l_{2}\right]$ with Dirichlet conditions at the boundaries,
    but with no additional condition at $x=0$.

    Conclude that this is the case in general. Namely for any graph, a
    vertex of degree two with Neumann vertex conditions ($\alpha=0$)
    is superfluous, i.e. one can erase the vertex and join the two incident
    edges to a one single edge such that the lengths add up.\\
    \textit{Hint}: You can read about this in \cite[section 2.2]{Berkolaiko_arXiv17}\\[0.1cm]
  \end{enumerate}
\end{itemize}

\textit{Comment}: Equation (\ref{eq:Secular_function_with_delta}) introduces a secular function which possess poles. The regularization process provides the secular function (\ref{eq:secular-function-with-delta-no-poles}) with no poles. This latter secular function is the more standard one. There is an explicit formula for it (with no need of regularization), which makes it also more convenient to use in proofs and analysis on quantum graphs. See more about that in Section \ref{sebsec:secular-func-second-approach}.

\subsection{From the interval to the star graph}

\textbf{~}\\[0.4cm]
\textbf{Background}: Consider a star graph which consists of one central
vertex which is connected to all other $V-1$ vertices.
We will enumerate the vertices as $i=0,1,\dots,E$ where $E=V-1$
is the total number of edges. The central vertex is privileged to
have the index $i=0$. The lengths of the edges are $\left\{ l_{01},\ldots,l_{0E}\right\} $.
The restrictions of a function $f$ to the edges are denoted by $\left\{ f_{01},\ldots,f_{0E}\right\} $.
The coordinates along the edges are chosen according to the convention
that $x_{0i}=0$ at the central vertex and $x_{0i}=l_{0i}$ at the
vertex $i$. We will assume Dirichlet boundary conditions at the vertices
$i=1,\dots,E$, and a repulsive vertex potential whose strength is
$\alpha_{0}\ge0$ at the central vertex. Namely, the vertex conditions
at the boundary vertices are $\delta$-type vertex conditions:
\[
\forall\,\,1\leq i\leq E\,\,;\,\,f_{0i}(l_{0i})=0,
\]
and the vertex conditions at the central vertex are
\[
\forall\,1\leq i<j\leq E\,\,;\,\,\,\,f_{0i}(0)=f_{0j}(0)\equiv f_{0}
\]
and
\[
\sum_{i=1}^{E}f_{0i}'(0)=\alpha_{0}f_{0}.
\]

Note that the problem considered in section \ref{subsec:Delta_type_conditions}
is a special case of the problem considered here with $E=2$.\\[0.2cm]

\noindent
\textbf{Questions:}\\[0.2cm]
\begin{itemize}
\item[\textbf{3.}] Show that a possible secular function for the graph is
  \begin{equation}
    \zeta(k)=\frac{\alpha_0}{k}+\sum_{i=1}^{E}\cot(kl_{0i}).\label{eq:secular_function_star}
  \end{equation}
  Namely, show that zeros of this function are all eigenvalues of the
  graph (note that there is not necessarily one to one correspondence
  - see the next question).

  Note that the secular function, (\ref{eq:Secular_function_with_delta})
  in the previous question is a special case of the one above (with
  $E=2$).\\[0.2cm]
\item[\textbf{4.}]
  What is the weakest assumption you need to assume in order to have
  a one to one correspondence between the zeros of the secular function
  (\ref{eq:secular_function_star}) and the graph's eigenvalues?\\[0.2cm]
\item[\textbf{5.}]
  Show that under the assumption you got in the previous question, the
  graph's spectrum is non degenerate. Namely, that each eigenvalue appears
  with multiplicity one.\\[0.2cm]
\item[\textbf{6.}] Under the same assumption as in questions \textbf{4} and \textbf{5}
  %(2) and (3),
  show the following interlacing
  properties:\\[0.2cm]
  \begin{enumerate}[(a)]
  \item The Dirichlet spectrum ($\alpha_{0}\rightarrow\infty$) of the star
    graph interlaces with the spectrum of the same star graph but with
    a Neumann vertex condition at the central point ($\alpha_{0}=0$).\\[0.1cm]
  \item For any positive value of $\alpha_{0}$ the $n$-the eigenvalue of
    the star graph is bounded from below by the $n$-th eigenvalue for
    $\alpha_{0}=0$ and from above by the $n$-th eigenvalue of the Dirichlet
    spectrum. \\[0.1cm]
  \item Similarly to the previous question, which interlacing property holds for the $n$-th eigenvalue of
    the star graph with a negative value of $\alpha_0$?\\[0.1cm]
  \end{enumerate}
  \textit{Comment}: See \cite[theorem 3.1.8]{BerKuc_graphs} for a general statement on eigenvalue interlacing for quantum graphs.  \\[0.2cm]

\item[\textbf{7.}]
  Assume now that $\alpha=0$ (Neumann vertex conditions at the center).
  Consider all star graphs with any number of edges $E$ and any positive edge
  lengths $\left\{ l_{01},\ldots,l_{0E}\right\} $, such that the total
  length of the edges, $L=\sum_{i=1}^{E}l_{0i}$ is fixed.\\[0.2cm]
  \begin{enumerate}[(a)]
  \item What is the supremum of the first eigenvalue among all choices of
    values for $E$ and $\left\{ l_{01},\ldots,l_{0E}\right\} $ (under
    the constraint above)? Is this supremum attained (i.e., is it a maximum)?\\[0.1cm]
  \item What is the infimum? Is it attained?\\[0.1cm]
  \end{enumerate}
  \textit{Hint}: The answer can be found in \cite{Fri_aif05}. \\
  \textit{Comment}: Similar questions of eigenvalue optimization on graphs and bounds on the eigenvalues are discussed in the recent works \cite{Ari_prep16, BanLev_ahp17, BerKenKurMug_jpa17, KenKurMalMug_ahp16}  \\[0.2cm]

\item[\textbf{8.}] \textit{This is a numerical exercise!} Choose any value for the number of edges,
  $E$, and any values for the edge lengths $\left\{ l_{01},\ldots,l_{0E}\right\} $.
  Plot the secular function for your choice and find its thirteen first
  zeros.\\[0.2cm]
\end{itemize}

\subsection{Secular function - first approach}\label{sebsec:secular-func-first-approach}

~\\[0.4cm]
\noindent
\textbf{Question:}\\[0.2cm]
\begin{itemize}
\item[\textbf{9.}]
  Consider an arbitrary quantum graph with $V$ vertices
  and $E$ edges. Assume there is an edge connecting vertices $1$ and $2$ and write the restriction of a function $f$ to this edge as
  \[
  f_{12}(x_{12})=A\cos kx_{12}+B\sin kx_{12}.
  \]
  Similarly, write
  \[
  f_{ij}(x_{ij})=\frac{f_{j}\sin(kx_{ij})+f_{i}\sin(k(l_{ij}-x_{ij}))}{\sin k l_{ij}},
  \]
  for the restriction of $f$ to an edge connecting vertices $i,\,j$ ($i<j$).
  Using this, obtain a set of homogeneous equations for the coefficients
  $A$, $B$, $f_{i}$ ($i=1,2,\dots V$) and derive a secular function
  that does not have poles at the Dirichlet spectrum of the edge $e=(1,2)$.

  \textit{Hint}: There is more than one solution to this - one may also
  reduce the number of equations (and variables) easily to $V$ (the
  number of vertices) without re-introducing poles.\\[0.2cm]
\end{itemize}

 \subsection{Secular function - scattering approach}\label{sebsec:secular-func-second-approach}

~\\[0.4cm]
\textbf{Background}: One may express the restriction of
the eigenfunction to the edge connecting vertices $i,j$ by
\begin{equation}
f_{ij}(x_{ij})=a_{ij}^{\mathrm{in}}e^{-ikx_{ij}}+a_{ij}^{\mathrm{out}}e^{ikx_{ij}}\ ,\label{eq:superposition_of_incoming_and_outgoing}
\end{equation}
where $a_{ij}^{\mathrm{in}}$ and $a_{ij}^{\mathrm{out}}$
are some coefficients. If the eigenfunction belongs to the eigenvalue
$k^{2}$, we can always choose values to those coefficients, such that
expression (\ref{eq:superposition_of_incoming_and_outgoing}) holds.

For a given vertex $i$ of degree $d_{i}$ we have $d_{i}$ coefficients
of the type $a_{ij}^{\mathrm{in}}$ and another $d_{i}$ coefficients
of the type $a_{ij}^{\mathrm{out}}$. Let us collect these into
$d_i$-dimensional vectors $\vec{a}^{(i),\mathrm{out}},~\vec{a}^{(i),\mathrm{in}}\in\mathbb{C}^{d_{i}}$.
The vertex conditions on
the vertex $i$ allow us to express the $\vec{a}^{(i),\mathrm{out}}$
as a linear transform of the $\vec{a}^{(i),\mathrm{in}}$
\begin{equation}
\vec{a}^{(i),\mathrm{out}}=\boldsymbol{\sigma}^{(i)}(k)\vec{a}^{(i),\mathrm{in}}. \label{eq:vertex-scattering-matrix-with-vectors}
\end{equation}
The components of the vectors
$\vec{a}^{(i),\mathrm{out}},~\vec{a}^{(i),\mathrm{in}}\in\mathbb{C}^{d_{i}}$
have the meaning
of incoming and outgoing wave amplitudes;
$\boldsymbol{\sigma}^{(i)}(k)$
is a unitary matrix of size $d_i \times d_i$
and is called the vertex scattering matrix.

To each quantum graph with $E$ edges one may associate
a unitary matrix $U(k)$ of dimension $2E \times 2E$,
known as the graph's \emph{quantum map}
or \emph{(discrete) quantum evolution operator}
that describes the connectivity of the graph, the matching
conditions at the vertices and the eigenvalue spectrum of the graph.
In Question~\textbf{12} the quantum map is derived explicitly
for a particular graph. The same procedure may be applied to
other quantum graphs to find the quantum map.\\[0.2cm]

%\pagebreak

\noindent
\textbf{Questions:}\\[0.2cm]
\begin{itemize}
\item[\textbf{10.}]  Show that the $\delta$-type vertex conditions (i.e., continuity of
  the eigenfunction at a vertex $i$ of degree $d_{i}$ and , $\alpha_{i}f_{i}=\sum_{j\in V_{i}}f_{ij}'(0)$)
  are equivalent to the vertex scattering matrix
  \begin{equation}
    {\sigma}_{jj'}^{(i)}(k)=\frac{2}{d_{i}+i\frac{\alpha_{i}}{k}}-\delta_{jj'}=\begin{cases}
      \frac{2}{d_{i}+i\frac{\alpha_{i}}{k}} & \quad j\neq j'\\
      \frac{2}{d_{i}+i\frac{\alpha_{i}}{k}}-1 & \quad j=j'
    \end{cases}.\label{eq:vertex-scattering-matrix}
  \end{equation}
  %\\[0.2cm]
  % \\
%\\
%\textbf{Question 2}:
\item[\textbf{11.}] Show that the vertex scattering matrix
  \[
  \boldsymbol{\sigma}^{(i)}=-\mathbf{1}_{d_{i}}+\frac{2}{v_{i}+i\frac{\alpha_{i}}{k}}\mathbb{E}_{d_{i}}
  \]
  is unitary. Here, $\mathbb{E}_{d_{i}}$ is the full $d_{i}\times d_{i}$
  matrix with all matrix elements equal to one. You may use $\mathbb{E}_{d_{i}}^{2}=d_{i}\mathbb{E}_{v_{i}}$
  and $\mathbb{E}^{*}=\mathbb{E}$ where $\mathbb{E}^{*}$ denotes the
  hermitian conjugate of $\mathbb{E}$.\\[0.2cm]
\item[\textbf{12.}]
%\textbf{Question 3}:
  Consider a star graph which consists of three
  edges. The central vertex is denoted by $0$ and supplied with Neumann
  vertex conditions. The boundary vertices are denoted by $1,2,3$ and
  are supplied with Neumann, Dirichlet, Dirichlet conditions, correspondingly.
  These notations and the edge lengths are shown in Fig. \ref{fig-star}.
  \begin{figure}[H]
    \includegraphics[scale=0.3]{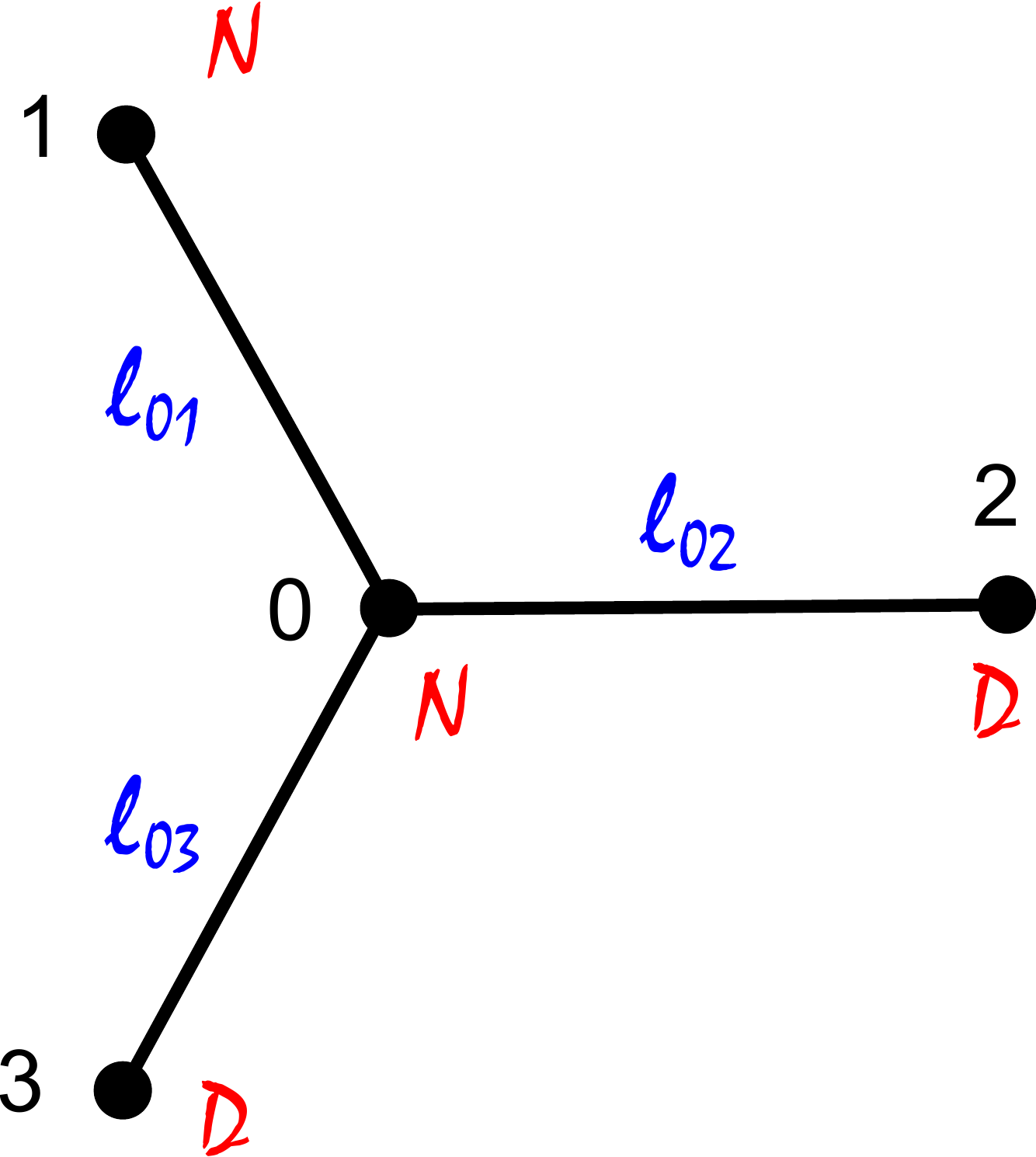}
    \caption{A star graph.}
    \label{fig-star}
  \end{figure}
  \noindent
  In this question you will explicitly build up the
  quantum map of the graph above by following the given
  sequence of instructions.\\[0.2cm]
  \begin{enumerate}[(a)]
  \item Write explicitly the scattering matrices, $\boldsymbol{\sigma}^{(i)}$,
    which correspond to each of the vertices $0,1,2,3$. You may compare
    with equation \eqref{eq:vertex-scattering-matrix} above.\\[0.1cm]
  \item Write the equation \eqref{eq:vertex-scattering-matrix-with-vectors}
    for each of the vertices. Use the explicit matrices which you have
    found in the previous section and write the components of the vectors
    $\vec{a}^{(i),\mathrm{in}}$ and $\vec{a}^{(i),\mathrm{out}}$ with
    explicit indices in each case (i.e., write $a_{0}^{(2),\mathrm{in}},\,
    a_{3}^{(0),\mathrm{out}}$,
    etc.).\\[0.1cm]
  \item Write (explicitly again) the ``big'' scattering matrix $S(k)$
    to fit the following set of equations
    \[
    \vec{a}^{\, \mathrm{out}}=S(k)\vec{a}^{\, \mathrm{in}},
    \]
    where
    \[
    \vec{a}^{\, \mathrm{in}}=\left(\begin{array}{c}
        a_{10}^{\mathrm{in}}\\
        a_{01}^{\mathrm{in}}\\
        a_{20}^{\mathrm{in}}\\
        a_{02}^{\mathrm{in}}\\
        a_{30}^{\mathrm{in}}\\
        a_{03}^{\mathrm{in}}
      \end{array}\right)\,\,\,\,\textrm{and}\,\,\,\,
    \vec{a}^{\, \mathrm{out}}=\left(\begin{array}{c}
        a_{01}^{\mathrm{out}}\\
        a_{10}^{\mathrm{out}}\\
        a_{02}^{\mathrm{out}}\\
        a_{20}^{\mathrm{out}}\\
        a_{03}^{\mathrm{out}}\\
        a_{30}^{\mathrm{out}}
      \end{array}\right).
    \]
    Remember that $S(k)$ merely consists of the different
    components of the single vertex scattering matrices $\boldsymbol{\sigma}^{(i)}(k)$
    and zero elements for edges that are not connected to each other.
    Pay special care to the order of the entries of the vectors above.\\[0.1cm]
  \item Write the matrix $T(k)$ such that it fits into the set
    of equations $$\vec{a}^{\, \mathrm{in}}=T(k)\vec{a}^{\, \mathrm{out}}$$ with $\vec{a}^{\, \mathrm{in}}$
    and $\vec{a}^{\, \mathrm{out}}$ as given above.\\
    \textit{Hint: you may use the fact $f_{ij}(x) = f_{ji}(l_{ji}-x)$}\\[0.1cm]
  \item A few tips to check yourself (no need to calculate, just in order
    to verify your answer).\\
    \textit{i.}
    The matrix $S(k)$ should be $k$-independent.\\
    \textit{ii.} The matrix $S(k)$should be unitary.\\
    \textit{iii.} The matrix $T(k)$ should be diagonal.\\[0.1cm]
  \end{enumerate}
  If you did all the above correctly, the quantum
  evolution operator is obtained by matrix multiplication of the two matrices,
  $U(k)=T(k)S(k)$.
  \\[0.2cm]
  % \textbf{Question 4}:
\item[\textbf{13.}]
  Consider a star graph with Neumann vertex conditions
  at the central vertex $i=0$ and Dirichlet vertex conditions at the
  boundary vertices $i=1,2,\dots,E$. Derive the quantum evolution map
  $U(k)$
  and show that the secular function can be reduced to the form
  \[
  \zeta(k)=
  \mathrm{det}\left(\mathbf{1}{}_{2E}-U(k)\right)=
  \mathrm{det}\left(\mathbf{1}{}_{E}
    +\widetilde{T}(2k)\boldsymbol{\sigma}^{(0)}\right),
  \]
  where $\boldsymbol{\sigma}^{(0)}$ is the central vertex scattering matrix and
  $\widetilde{T}(k)$ is a diagonal $E\times E$ matrix,
  $T(k)_{ee'}=\delta_{ee'}e^{ikl_{e}}$. You can gain a good intuition
  for the solution of this question from your solution to the previous
  question.\\[0.2cm]
  % \\
  % \\
  % \textbf{Question 5}:
\item[\textbf{14.}]
  Show that the following secular function is real
  \begin{equation}
  \tilde{\zeta}(k)=\sqrt{\mathrm{det}\left(S^{*}(k)T^{*}(k)\right)}\mathrm{det}\left(\mathbf{1}_{2E}-U(k)\right). \label{eq:real-secular-function}
  \end{equation}
  Remember that $U\left(k\right)=T(k)S(k)$ and use the unitarity of
  $T(k)$ and $S(k)$.\\
  \textit{Comment}: The secular function $\tilde{\zeta}(k)$ above may
  even be differentiable in $k$, if the complex branch of the square
  root is appropriately chosen.\\[0.2cm]
\end{itemize}

%~

\section{Trace Formula and Periodic Orbits}

\subsection{The Trace Formula for the Spectrum of a Unitary Matrix\vspace{4mm}}

\textbf{~}\\
\textbf{Background:} Consider an $M\times M$ unitary matrix $U$ with
unimodular eigenvalues $e^{i\theta_{\ell}}$ for $\ell=1,\dots,M$.
 One may extend the spectrum of eigenphases $\theta_{\ell}$ periodically
beyond the interval $0\le\theta<2\pi$. The extended spectrum then
consists of the numbers
\[
\theta_{\ell,n}=\theta_{\ell}+n 2\pi\qquad n\in\mathbb{Z}\ .
\]
Assume that $\theta_{\ell}\neq0$ and $\theta\neq\theta_{\ell,n}$
is real. \\[0.2cm]

\noindent
\textbf{Question:}\\[0.2cm]
\begin{itemize}
\item[\textbf{15.}]
  Consider the spectral
  counting function
  \begin{align*}
    N(\theta)= & \sum_{n=0}^{\infty}\sum_{\ell=1}^{M}\vartheta(\theta-\theta_{\ell}-n2\pi),
  \end{align*}
  where $\vartheta$ is the Heaviside step function.
  Show that one may write it as the following trace formula.
  \[
  N(\theta)=\frac{M\theta}{2\pi}-\frac{1}{\pi}\mathrm{Im}\ \mathrm{log}\ \mathrm{det}\ (1-U)+\frac{1}{\pi}\mathrm{Im}\ \mathrm{log}\ \mathrm{det}\ (1-e^{-i\theta}U).
  \]

  \textit{Comment:} The definition of the spectral counting function
  may be extended to $\theta=\theta_{n,\ell}$ such that both expressions
  (the defining expression and the trace formula) remain consistent.
  For this one replaces the last term in the trace formula by the limit
  \[
  \frac{1}{\pi}\mathrm{Im}\ \mathrm{log}\ \mathrm{det}\ (1-e^{-i\theta}U)\mapsto\lim_{\epsilon\to0^{+}}\frac{1}{\pi}\mathrm{Im}\ \mathrm{log}\ \mathrm{det}\ (1-e^{-i\theta-\epsilon}U)
  \]
  and sets $\vartheta(0)=1/2$. Replacing $U\mapsto e^{-\epsilon}U$
  and considering the limit $\epsilon\to0^{+}$ also helps to regularize
  certain expansions that may come up in the proof of the trace formula
  because the trace formula for finite $\epsilon>0$ does not have any
  singularities for $\theta$ on the real line.

  We strongly recommend plotting the regularized expression for $N(k)$
  with a (small) positive value for $\epsilon$ for a given
  unitary matrix $U$ (which may be chosen diagonal).

  \textit{Hint}: There are several ways to perform the derivation.
  One interesting derivation is based on Poisson summation. This method
  requires $\epsilon$-regularization as mentioned in the comment above.\\[0.2cm]
  \begin{enumerate}[(a)]
  \item Write the spectral counting function as
    \begin{align*}
      N(\theta)= & \sum_{n=-\infty}^{\infty}\sum_{\ell=1}^{M}\vartheta(\theta-\theta_{\ell}-n2\pi)\vartheta(\theta_{\ell}+n2\pi).
    \end{align*}
  \item The Poisson summation formula for a smooth function $f(x)$ which
    decays sufficiently fast for $|x|\to\infty$ (so that all sums and
    integrals converge absolutely) reads
    \[
    \sum_{n=-\infty}^{\infty}f(n)=\sum_{\nu=-\infty}^{\infty}\int_{-\infty}^{\infty}e^{2\pi i\nu x}f(x)dx\ .
    \]
    We want to apply it to $f(x)=\sum_{\ell=1}^{M}\vartheta(\theta-\theta_{\ell}-2\pi x)\ \vartheta(\theta_{\ell}+2\pi x)$
    which is not smooth. In this case the Poisson sum is not absolutely
    convergent but may be regularized by introducing an additional factor
    $e^{-\epsilon|\nu|}$ and taking the limit $\epsilon\to0^{+}$.\\[0.1cm]
  \item While using the formula above, evaluate separately the $\nu=0$ term
    from the other terms.\\[0.1cm]
  \item Compare this term by term to the expansion
    \[
    \mathrm{log}\ \mathrm{det}\ (1-e^{-i\theta}U)=\mathrm{tr}\ \mathrm{log}\ (1-e^{-i\theta}U)=-\lim_{\epsilon\to0^{+}}\sum_{n=1}^{\infty}\frac{1}{n}e^{-in\theta-n\epsilon}\mathrm{tr}\ U^{n}
    \]
    (or the complex conjugate version) in order to perform the sum over
    $\nu$ in the Poisson summation.\\[0.1cm]
  \end{enumerate}

  \textit{Comment:} You may compare this to the trace formula of a quantum
  graph. Do this by showing that the above formula is equivalent to
  the one obtained in \cite{GnuSmi_ap06} for quantum graphs if all
  the edge lengths of the graph are the same. Note that the derivation
  of the trace formula in this reference uses a different (more general)
  method.\\[0.2cm]
\end{itemize}

\subsection{Periodic Orbits}

~\\[0.4cm]
\noindent
\textbf{Background:} Consider a quantum graph with $E$ edges of lengths
$l_{e}$ ($e=1,\dots,E$) with Neumann matching conditions and let
$U(k)$ be the unitary, $k$-dependent $2E\times2E$ matrix
\[
U(k)_{\alpha'\alpha}=e^{ikl_{\alpha'}}S_{\alpha'\alpha}
\]
where $S_{\alpha'\alpha}$ is the scattering amplitude from the directed
edge $\alpha$ to directed edge $\alpha'$.

Note that $S_{\alpha'\alpha}=0$
unless the end vertex of the directed edge $\alpha$ coincides with
the start of $\alpha'$ \textendash{} one then says that $\alpha'$
follows $\alpha$.

We have
\[
\mathrm{tr}\ U(k)^{n}=\sum_{\alpha_{1},\dots,\alpha_{n}=1}^{2E}e^{ikl_{\alpha_{1}}}S_{\alpha_{1}\alpha_{n}}e^{ikl_{\alpha_{n}}}S_{\alpha_{n}\alpha_{n-1}}\dots e^{ikl_{\alpha_{3}}}S_{\alpha_{3}\alpha_{2}}e^{ikl_{\alpha_{2}}}S_{\alpha_{2}\alpha_{1}}\ ,
\]
where $\left\{ \alpha_{1},\dots,\alpha_{n}\right\} $ is a set of
$n$ directed edges of the graph.
In the sum above, each term is of the form $A_{\gamma}e^{ik\ell_{\gamma}}$
\begin{align*}
  A_{\gamma}&=  S_{\alpha_{1}\alpha_{n}}S_{\alpha_{n}\alpha_{n-1}}\dots S_{\alpha_{2}\alpha_{1}}\\
  \ell_{\gamma}&=  l_{\alpha_{n}}+l_{\alpha_{n-1}}+\dots l_{1}
\end{align*}
Here $\gamma=\left(\alpha_{1},\dots,\alpha_{n}\right)$ is a fixed
set of summation indices which corresponds to a sequence of directed
edges.

\begin{figure}[H]
\hfill{}\includegraphics[scale=1.1]{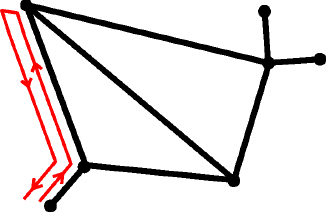}\hfill{}\includegraphics[scale=1.1]{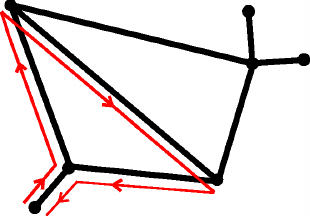}\hfill{}\includegraphics[scale=1.1]{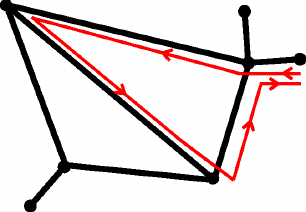}\hfill{}\caption{Three examples of closed trajectories on a graph.}
\label{fig:periodic_orbits}
\end{figure}

Note that $A_{\gamma}\neq 0$ only if the edge $\alpha_{j+1}$ follows
the edge $\alpha_{j}$ in the graph, for $j=1,\dots,n$ (in this context
$\alpha_{n+1}\equiv\alpha_{1}$) \textendash{} i.e. only if $\gamma$
is a closed trajectory on the graph. For such a trajectory, $l_{\gamma}$
is its total length, i.e., the sum of its edge lengths. By definition,
the \textbf{\textit{length spectrum}} of the graph is the set of all
lengths $\{l_{\gamma}\}$ of closed trajectories.
 You can see a few examples of such
closed trajectories in Fig.~\ref{fig:periodic_orbits}.

We now list a few observations and definitions related to periodic
orbits.\\[0.2cm]
\begin{itemize}
\item Note that any cyclic permutation, $\gamma'$, of the indices in the
  trajectory $\gamma$ gives a different closed trajectory with the
  same contribution $A_{\gamma'}=A_{\gamma}$ and $\ell_{\gamma'}=\ell_{\gamma}$.\\[0.1cm]
\item The equivalence class $\overline{\gamma}=\overline{\alpha_{1}\dots\alpha_{n}}$
  that contains all cyclic permutations of a given closed trajectory
  $\gamma=(\alpha_{1},\dots,\alpha_{n})$ is called a \textbf{\textit{periodic
      orbit}} with \textbf{\textit{period}}\textbf{ $n$} on the graph.\\[0.1cm]
\item The periodic orbit $\overline{\gamma}=\overline{\alpha_{1}\dots\alpha_{n}}$
  is a \textit{ }\textbf{\textit{primitive periodic orbit}}\textbf{
  }of primitive period $n$ if the sequence of indices $\alpha_{1},\dots,\alpha_{n}$
  is not a repetition of a shorter sequence.\\[0.03cm] All the closed
  trajectories in Fig. \ref{fig:periodic_orbits} represent primitive
  periodic orbits.\\[0.1cm]
\item If $\overline{\gamma}$ is a periodic orbit with period $n$, then
  there exists a unique primitive periodic orbit $\overline{\gamma}_{p}$
  with \textbf{\textit{primitive period}} $n_{p}$ such that $\overline{\gamma}$
  is a repetition of $\overline{\gamma}_{p}$ and $n=rn_{p}$. Here
  $r\ge1$ is the integer \textbf{\textit{repetition number}} of $\overline{\gamma}$
  and $n_{p}$ the \textbf{\textit{primitive period}} of $\overline{\gamma}$.\\[0.1cm]
\item If $r=1$ then $\overline{\gamma}=\overline{\gamma}_{p}$ and $\overline{\gamma}$
  is primitive.\\[0.1cm]
\item We write $\overline{\gamma}=\overline{\gamma}_{p}^{\ r}$ for the
  $r$-th repetition of the primitive orbit $\overline{\gamma}_{p}$.\\[0.1cm]
\end{itemize}

\textit{Comment}: The trace formula of the spectral counting function may be expressed in terms of an infinite sum over the graph's periodic orbits.
To answer the following questions, you should read more about that, for example in \cite[section 3.7.4]{BerKuc_graphs} or \cite[section 5.2]{GnuSmi_ap06}. \\[0.2cm]
% Further reading is available from the review paper ``Quantum graphs:
% Applications to quantum chaos and universal spectral statistics''
% by S. Gnutzmann and U. Smilansky (see link in the course website).

\noindent
\textbf{Questions:}\\[0.2cm]
\begin{itemize}
\item[\textbf{16.}]
  Consider the dihedral graphs as given in Fig.~\ref{dihedral}
  \begin{figure}[H]
    \hfill{}\includegraphics[scale=0.6]{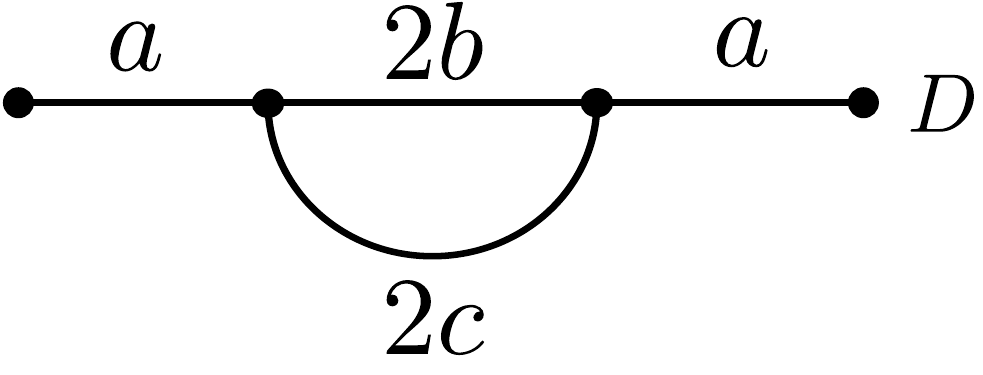}\hfill{}\includegraphics[scale=0.6]{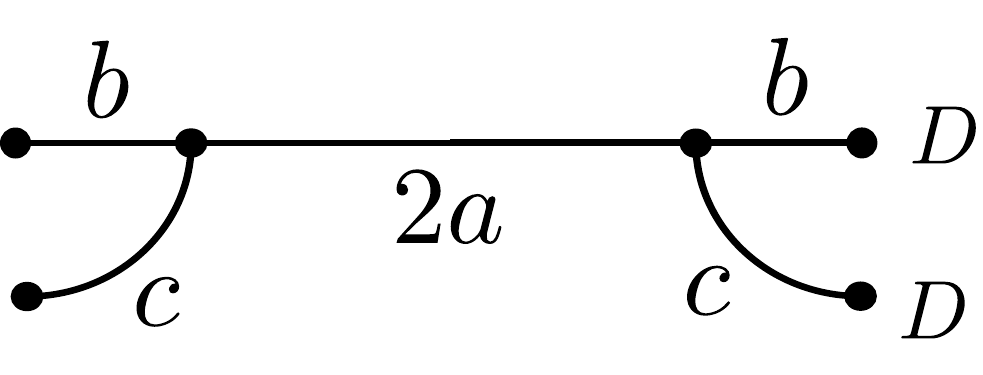}\hfill{}\caption{Two isospectral quantum graphs. Their edge lengths are indicated by
      the parameters $a,b,c$. The vertices marked with `D' have Dirichlet
      conditions and all other vertices have Neumann conditions. \label{dihedral}}
  \end{figure}
  \noindent It is known that those graphs are isospectral \cite{BanShaSmi_jpa06}.
  Solve the following `paradox':\\
  The isospectrality of those graphs means that their spectral counting functions are equal. Hence, the periodic orbit expansions of those counting functions are the same. Therefore, both graphs should have the same set of periodic orbits. Nevertheless, the graph on the left
  has a periodic orbit of length $2a$, whereas the graph on the right
  does not have such an orbit. How is this possible?\\[0.2cm]
\item[\textbf{17.}]
  %\noindent \textbf{Question 2}
  This question concerns the tetrahedron graph, i.e.
  the complete Neumann graph with $V=4$ vertices.
  \begin{enumerate}[(a)]
  \item \emph{Warm-up:} Choose two periodic
  orbits on the graph, such that one of them is primitive and the second
  is some repetition of the first. For each of those periodic orbits,
  $\gamma$, evaluate the following quantities:\\[0.1cm]
    \begin{enumerate}[(i)]
    \item The period $n$ of the orbit.\\[0.1cm]
    \item The length $\ell_{\gamma}$ of the orbit (expressed in terms of the graph
        edge lengths).\\[0.1cm]
    \item The coefficient $A_{\gamma}$ which corresponds to the orbit in the periodic
        orbits expansion (write the explicit number).\\[0.1cm]
    \item The primitive period $n_{p}$ and the repetition number $r$ ?\\[0.1cm]
    \end{enumerate}
  \item Assume that all bond lengths are incommensurate. Go over all periodic orbits of period $n=5$ and write their contribution to the length
    spectrum $\sigma_{\ell}$ (expressed in terms of the graph
        edge lengths).
    What are the corresponding quantum amplitudes $A_{\gamma}$ of those orbits?\\[0.1cm]
  \item Now assume that all edge lengths are equal.
    How does your answer to the previous question change?\\[0.1cm]

  \item The connectivity matrix $C$ of a simple graph with
    $V$ vertices is the real symmetric  $V\times V$ matrix with
    entries
    \[
  C_{i,j}=\begin{cases}
    1 & \text{vertices $i$ and $j$ are connected by an edge,}\\
    0 & \text{else.}
  \end{cases}
  \]
  Note that for simple graphs there are
    no loops, i.e. edges that connect a vertex to itself, so $\forall i, \, C_{ii}=0$.
    The connectivity matrix may be used to count the number of trajectories
    that connect the vertices $i$ and $j$ in $n$ steps \textendash{}
    this number is just $\left[C^{n}\right]_{ji}$. \\
    Consider the connectivity matrix of the tetrahedron, $C_{ij}=1-\delta_{ij}$. \\
    Show that as $n\to\infty$, $\left[C^{n}\right]_{ji}\sim ce^{\alpha n}$ and find $\alpha$.

    \textit{Hint}: diagonalize $C$. \\[0.1cm]
  \end{enumerate}

  \textit{Comment}:  The last question shows that there is an exponential growth in the number of orbits on the graph. This holds in particular for periodic orbits and makes it difficult to use the periodic orbit expansion of the trace formula for spectral computations.
  \\[0.2cm]

  %
%\textbf{Question 3:}
\item[\textbf{18.}]
  In the previous question we counted the number
  of trajectories between two vertices via the (vertex) connectivity
  matrix. An alternative approach is based on the $2E\times2E$ edge
  adjacency matrix $\mathcal{B}$ whose indices correspond to the directed
  edges where
  \[
  \mathcal{B}_{\alpha\alpha'.}=\begin{cases}
    1 & \text{if \ensuremath{\alpha}\, follows \ensuremath{\alpha'},}\\
    0 & \text{else.}
  \end{cases}
  \]
  The main difference is that $\left[\mathcal{B}^{n}\right]_{\alpha\alpha'}$
  counts the number of trajectories that start on the directed edge
  $\alpha'$ and end after $n$ steps on the directed edge $\alpha$.
  Both approaches can also be used to count the number of periodic orbits
  via the traces. We will explore this here for the edge connectivity
  matrix $\mathcal{B}$.\\[0.2cm]
  \begin{enumerate}[(a)]
  \item Show that $\frac{1}{n}\mathrm{tr}\thinspace\mathcal{B}^{n}=\sum_{\gamma:n_{\gamma}=n}\frac{1}{r_{\gamma}}$
    where the sum is over all periodic orbits of period $n$. Conclude
    that if $n$ is a prime number then $\frac{1}{n}\mathrm{tr}\thinspace\mathcal{B}^{n}$
    is the number of periodic orbits of period $n$.\\[0.1cm]
  \item Derive an expression for the number of periodic orbits of period $n$
    in terms of traces of powers of $\mathcal{B}$ for \\
    (\textit{i.}) $n=p^{j}$
    where $p$ is a prime number and $j\ge2$ an integer, and \\
    (\textit{ii.})
    $n=p_{1}p_{2}$ where $p_{1}$ and $p_{2}$ are prime numbers.\\
    Make an educated guess for the general expression when $n$ has the prime
    number decomposition $n=\prod_{m}p_{m}^{j_{m}}$ where $j_{m}\ge0$
    is the multiplicity of the $m$-th prime. \\[0.2cm]
  \end{enumerate}
  % ~
  %
  % \textbf{Question 4.}
\item[\textbf{19.}]
  This question demonstrates that given some lengths
  of periodic orbits of an unknown graph, one can reconstruct the graph.
  \\[0.2cm]
  \begin{enumerate}[(a)]
  \item Find the graph with the following properties: \\
    (\textit{i.}) the total length (sum of all edge lengths) is $1\frac{5}{6}$;\\
    (\textit{ii.}) the lengths of all periodic orbits whose length is not
    greater than $5$ are given by the list
    \[
    \frac{2}{3},\,1\frac{1}{3},\,2,\,2\frac{2}{3},\,3,\,3\frac{1}{3},\,3\frac{2}{3},\,4,\,4\frac{1}{3},\,4\frac{2}{3},\,5.
    \]
    Draw the graph and indicate the edge lengths on the drawing.\\[0.1cm]
  \item Find the graph with the following properties: \\
    (\textit{i.}) the total length is $5\frac{13}{15}$;\\
    (\textit{ii.}) the lengths of all periodic orbits whose length is not
    greater than $5$ are given by the list
    \[
    2,\,2\frac{1}{3},\,2\frac{1}{6},\,2\frac{2}{5},\,4,\,4\frac{1}{6},\,4\frac{1}{3},\,4\frac{1}{3},\,4\frac{2}{5},\,4\frac{1}{2},\,4\frac{17}{30},\,4\frac{2}{3},\,4\frac{11}{15},\,4\frac{4}{5},\,5.
    \]
    Note that there are two \uline{different} periodic orbits of length $4\frac{1}{3}$. Also, any number which appears only once in the list above indicates that there is exactly one periodic orbit of that length. \\
    What is the graph this time?\\[0.1cm]
  \item Try to think how to construct a general algorithm for finding the
    graph out of knowing its total metric length and lengths of all of
    its periodic orbits.\\
    Assume that the graph is simple (no loops and no multiple edges)
    and that its edge lengths are incommensurate.\\[0.1cm]
    % \item \RB{I think that Kurasov did something by removing one of the restrictions above. Maybe this can also be given as an exercise}
  \end{enumerate}
  \noindent \textit{Hint:} What is the shortest periodic orbit of a
  graph?\\
  \textit{Another hint:} The answer can be found in \cite{GutSmi_jpa01}.\\[0.2cm]
\end{itemize}

\subsection{The constant term of the Trace formula}

\textbf{~}\\[0.4cm]
\textbf{Background:} The spectral counting function of a quantum graph is
\[
N(k):=\left| \left\{ \lambda\in\mathbb{R} \textrm{ is an eigenvalue}  \, : \, \lambda < k^2  \right\} \right|,
\]
where eigenvalues are counted with their multiplicity.
One of the forms of the trace formula for the spectral counting function is
\begin{equation}
    N(k)=
    N_0+\frac{\mathcal{L}}{\pi}k-\mathrm{\lim_{\epsilon\rightarrow0}\frac{1}{\pi}Im}\ \mathrm{log}\ \tilde{\zeta}(k+i\epsilon), \label{eq:trace-formula-w-const-term}
\end{equation}
 where $\tilde{\zeta}(k)$ is the real secular function given in (\ref{eq:real-secular-function}), and $N_0$ is a constant term.
 The last two terms in (\ref{eq:trace-formula-w-const-term}) equal to the number of real zeros of $\tilde{\zeta}$ with absolute value smaller than $k$. Those zeros are in one to one correspondence with the graph eigenvalues (including multiplicity), with the exception of $k=0$. The value of $\tilde{\zeta}$ at $k=0$ does not correspond to the multiplicity of the zero eigenvalue and this `mismatch' is compensated by the constant term $N_0$ in (\ref{eq:trace-formula-w-const-term}). \\
 The expression for this term was originally derived in \cite[lemma 1]{KurNow_jpa05},\cite{KurNowCor_jpa06}. Other works related to this subject are \cite{FulKucWil_jpa07, HarrWey_lmp18}. Further reading in \cite[section 3.7]{BerKuc_graphs} and \cite[section 5]{GnuSmi_ap06} is recommended.\\[0.2cm]

\noindent
\textbf{Question:}\\[0.2cm]
\begin{itemize}
\item[\textbf{20.}]
  In the following question we derive the value of $N_0$ for a Neumann graph (all vertex matching conditions are of Neumann
aka Kirchhoff type). Initially, assume that the graph has a
single connected component. Let $\vec{E}$ be the space of directed edges
  on the graph (this space is of dimension $2E$, where $E$ is the
  number of edges).\\[0.2cm]
  \begin{enumerate}[(a)]
  \item Let $\omega:\vec{E}\rightarrow\mathbb{C}$ such that
    \[
    \forall\left(i,j\right)\in\vec{E}\quad \omega\left(i,j\right)=-\omega\left(j,i\right)
    \]
    and
    \[
    \forall i\quad\sum_{j\sim i}\omega\left(i,j\right)=0,
    \]
    where $j\sim i$ means that the vertex $j$ is adjacent to the vertex
    $i$. All such functions $\omega:\vec{E}\rightarrow\mathbb{C}$ form
    a vector space (over $\mathbb{C}$). Prove that the dimension of this
    space is $\beta:=E-V+1$.

    \textit{Hint}: Start by considering a tree graph. \\[0.1cm]
  \item Let $\vec{a}^{in}\in\mathbb{C}^{2E}$ with entries denoted by $a_{j}^{\left(i\right),in}$
    (for $i\sim j$), such that the following is satisfied
    \[
    \forall i\sim j\,,\,i\sim k\,\,\quad a_{j}^{\left(i\right),in}+a_{i}^{\left(j\right),in}=a_{k}^{\left(i\right),in}+a_{i}^{\left(k\right),in}
    \]
    and
    \[
    \forall i\quad\sum_{j\sim i}\left(-a_{j}^{\left(i\right),in}+a_{i}^{\left(j\right),in}\right)=0.
    \]
    Prove that the dimension of the vector space which contains all such
    solutions $\vec{a}^{in}\in\mathbb{C}^{2E}$ is $\beta+1$.\\[0.1cm]
  \item Note that you have shown $\dim\ker\left(\mathbf{1}-S\right)=\beta+1$, which implies $N_{0}=\frac{1-\beta}{2}$.
    Show that for a Neumann graph with $C$ (disjoint) connected components, the constant term of the trace formula (\ref{eq:trace-formula-w-const-term}) is $N_0=\frac{C-\beta}{2}$. \\
    You can use the generalized definition of $\beta$, which is $\beta:=E-V+C$
    (this value can be obtained by summing over all the $\beta$'s of
    the different components). \\[0.2cm]
  \end{enumerate}
\end{itemize}

\section{Further Topics}

\subsection{Quantum to Classical correspondence for Quantum Graphs}

~\\[0.4cm]
\noindent
\textbf{Background:} In this question we study the classical dynamics
of a quantum graph. This will help us to understand in what sense
the classical dynamics that corresponds to a quantum graph is 'chaotic'.
Remember that the quantum evolution map, $U(k)$, contains amplitudes
for scattering processes to go from one directed edge to another.
We define a corresponding classical map, $M$, by replacing the amplitudes
$U(k)_{\alpha\alpha'}$ by
\[
M_{\alpha\alpha'}=\left|U(k)_{\alpha\alpha'}\right|^{2}=\left|S_{\alpha\alpha'}\right|^{2}.
\]
Hence, $M$ is a matrix of dimensions $2E\times2E$, which contains the probabilities
for the scattering events.\\[0.2cm]
\noindent
\textbf{Question:}\\[0.2cm]
\begin{itemize}
\item[\textbf{21.}]
  By following the steps below prove that the
  matrix $M$ defines a Markov process on the set of directed edges
  with the stated additional properties.\\[0.2cm]
  \begin{enumerate}[(a)]
  \item Prove that the matrix $M$ is a \textit{bi-stochastic (doubly stochastic)
      matrix}\textit{\emph{. Namely, prove that}}
    \[
    \sum_{\alpha=1}^{2B}M_{\alpha\alpha'}=\sum_{\alpha'=1}^{2B}M_{\alpha\alpha'}=1.
    \]
  \item Use the bi-stochastic property to verify that the following definition
    of a Markov process on the directed edges of the graph is well-defined.
    Let $P_{\alpha}(n)$ be the probability to find a particle on the directed
    edge $\alpha$, at some (discrete) time $n$. We can then define the
    probabilities to find the particle on the directed edge $\alpha$,
    at time $n+1$, by
    \[
    P_{\alpha}(n+1)=\sum_{\alpha'}M_{\alpha\alpha'}P_{\alpha'}(n)
    \]
    or, in short $P(n+1)=MP(n)$. In particular show that if $P(n)$ satisfies
    $\sum_{\alpha}P_{\alpha}(n)=1$ and $P(n)_{\alpha}\ge0$, then $P(n+1)$
    satisfies the same properties. That is if $P(n)$ is a probability
    vector then $P(n+1)$ is a probability vector.\\[0.2cm]
  \end{enumerate}
\end{itemize}

\noindent
\textbf{Further background:} We next consider the equilibration properties
of the Markov process $P(n+1)=MP(n)$. Let $P^{\mathrm{inv}}=\frac{1}{2E}$
be the equi-distributed probability vector on the graph. For any quantum
graph this is an invariant probability vector, i.e. $$M\,P^{\mathrm{inv}}=P^{\mathrm{inv}}$$

The classical dynamics which corresponds to a quantum graph is chaotic
in the following sense: The Markov process on the graph is called
\textit{ergodic} if
\[
\lim_{n\rightarrow\infty}\frac{1}{n}\sum_{m=1}^{n}P(m)=P^{\mathrm{inv}}
\]
for every initial probability vector $P(0)$.
We call a graph dynamically connected\footnote{In the literature a matrix $M$ with this property is sometimes called irreducible.}, if for any two directed edges
$\alpha$ and $\alpha'$ there is an integer $n>0$ such that
$ (M^n)_{\alpha \alpha'} \neq 0$ (i.e one can get from one directed edge to
another with a non-vanishing probability in a finite number of steps).
Every dynamically connected graph
is ergodic.
Most graphs are also\textit{ mixing} which is the stronger property
\[
\lim_{n\rightarrow\infty}P(n)=P^{\mathrm{inv}}
\]
for every initial probability vector $P(0)$.\\[0.2cm]

\noindent
\textbf{Question:}\\[0.2cm]
\begin{itemize}
\item[\textbf{22.}]
  In this question we will characterize ergodicity and mixing
  in terms of the eigenvalue spectrum of the bi-stochastic matrix $M$.
  The results apply to any bi-stochastic Markov process on a directed
  graph whether or not there is a corresponding quantum graph such that
  $M_{\alpha\alpha'}=\left|U(k)_{\alpha\alpha'}\right|^{2}$ in terms of the
  quantum map.\\[0.2cm]
  \begin{enumerate}[(a)]
  \item Prove that all the eigenvalues of $M$ are either on the unit circle
    or inside it. Namely, if we denote the set of eigenvalues by $\left\{ \lambda_{i}\right\} $,
    then $\forall i\,;\,\,\left|\lambda_{i}\right|\leq1$.\\[0.1cm]
  \item We know that $M$ has at least one eigenvalue which equals 1 (the
    corresponding eigenvector is $P^{\mathrm{inv}}$). Let us denote this
    eigenvalue by $\lambda_{1}$ ($\lambda_{1}=1$). Prove that if $\min_{2\leq i\leq2E}\left(1-\left|\lambda_{i}\right|\right)>0$
    then the graph is mixing.

    \textit{Hint}: It might be useful to prove the convergence property using
    the vectors $L_{1}$-norm.\\[0.1cm]
  \item Using the notation above ($\lambda_{1}=1$), prove that if $\min_{2\leq i\leq2E}\left(\left|1-\lambda_{i}\right|\right)>0$
    then the graph is ergodic.

    \textit{Hint}: Again, use the $L_{1}$-norm .\\[0.1cm]
  \end{enumerate}

  \textit{Comment:} Note that the conditions above
  are consistent with the trivial fact
  that mixing is a stronger notion than ergodicity (namely, that every
  mixing system is also ergodic which follows directly from the
  definition).\\
  The quantity $\Delta:=\min_{2\leq i\leq2E}\left(\left|1-\lambda_{i}\right|\right)$
  is called the spectral gap and it determines the convergence rate
  of $\lim_{n\rightarrow\infty}\frac{1}{n}\sum_{m=1}^{n}P(m)$ (the
  greater the gap, the quicker is the convergence).\\
  Similarly, $\tilde{\Delta}:=\min_{2\leq i\leq2E}\left(1-\left|\lambda_{i}\right|\right)$
  determines the convergence rate of $\lim_{n\rightarrow\infty}P(n)$.\\[0.2cm]
\end{itemize}

\subsection{{The quadratic form}}

~\\[0.4cm]
\noindent
\textbf{Background:} Consider a quantum graph with the edge and vertex sets, $\mathcal{E}$ and $\mathcal{V}$.
We take the operator to be the Laplacian with $\delta$-type
  vertex conditions (see (\ref{eq:Delta_type_conditions_general})).
The quadratic form of this operator is:
\[
h[f] := \sum_{e\in \mathcal{E}} \int_{0}^{l_e} \frac{{\rm d}f}{{\rm d}x_{e}} \overline{\frac{{\rm d}g}{{\rm d}x_{e}}}{\rm d}x_{e}+ \sum_{v\in\mathcal{V}\,;\,\alpha_v\neq\infty} \alpha_v f(v) \overline{g(v)},
\]
where $\alpha_v$ is the coupling coefficient of the $\delta$-type vertex condition at vertex $v$, and $\alpha_v=\infty$ indicates Dirichlet vertex conditions. The length of the edge $e$ is denoted by $l_e$.\\
The domain $D(h)$ of this quadratic form consists of all functions $f$ on the metric graph that satisfy the following three conditions:\\
(\textit{i.}) for each edge $e$ the restriction $f|_e$ belongs to the
Sobolov space
$H^1([0,l_e])$, \\
(\textit{ii.}) $f$ is continuous at each vertex, and \\
(\textit{iii.}) $f(v)=0$ at each vertex $v$ for which $\alpha_v=\infty$.

The quadratic form is useful for variational characterization of the spectrum.
More on this topic is found in \cite[section 1.4.3]{BerKuc_graphs}.
\\[0.2cm]

\noindent
\textbf{Questions:}\\[0.2cm]
\begin{itemize}
\item[\textbf{23.}]
  Consider a Schr\"odinger operator with
  a bounded non-negative potential ($V\geq0$) on a quantum graph, $H\psi=-\psi''+V\psi$.
  Show that the spectrum of this operator is non-negative if all the
  vertex conditions are of $\delta$-type with non-negative coupling
  coefficients (i.e., $\forall\,v\in\mathcal{V}\,\,;\,\,\,\alpha_{v}\geq0$).\\
  \textit{Hint}: You need to modify the quadratic form given above to fit the case of Schr\"odinger operator with a potential. \\[0.2cm]
\item[\textbf{24.}]
  Prove the following statements:
\begin{enumerate}[(a)]
  \item Let $\lambda=\lambda(\alpha)$ be a simple eigenvalue of
  a graph with $\delta$-type vertex condition
  at a certain vertex $v$ with the coupling coefficient $\alpha\neq\infty$.
  The operator on this graph is just the Laplacian (no potential). Then
  \[
  \frac{{\rm d}\lambda}{{\rm d}\alpha}=\left|f\left(v\right)\right|^{2}.
  \]
  \item Now, re-parameterize the  $\delta$-type vertex condition at $v$ as:
  \[
  \zeta\sum_{e\in E_{v}}\frac{{\rm d}f}{{\rm d}x_{e}}(v)=-f(v),
  \]
  with $E_{v}$ denoting the set of edges adjacent to $v$.
  This parametrization allows Dirichlet condition ($\zeta=0$) and excludes Neumann condition ($\zeta=\infty$).
  Show that if the simple eigenvalue is now given by $\lambda=\lambda(\zeta)$ then the derivative is
  \[
  \frac{{\rm d}\lambda}{{\rm d}\zeta}=\left|\sum_{e\in E_{v}}\frac{{\rm d}f}{{\rm d}x_{e}}\left(v\right)\right|^{2}.
  \]
    \end{enumerate}~\textit{Hint}: The answer is given in \cite[proposition 3.1.6]{BerKuc_graphs} \\[0.2cm]
\end{itemize}

%\pagebreak
\subsection{{From quantum graphs to discrete graphs}}

~\\[0.4cm]
\noindent
\textbf{Question:}\\[0.2cm]
\begin{itemize}
\item[\textbf{25.}]
  In this question we consider the spectral connection between quantum graphs
  and discrete graphs.\\[0.2cm]
  \begin{enumerate}[(a)]
  \item Consider an arbitrary quantum graph with $V$ vertices and $E$ edges.
    Assume that Neumann conditions are imposed at all vertices and that
    all edges are of the same length, $l$. Use the following representation
    for an eigenfunction with eigenvalue $k^{2}$ on the edge $\left(i,j\right)$
    \[
    f_{ij}(x_{ij})=\frac{f_{j}\sin(kx_{ij})+f_{i}\sin(k(l-x_{ij}))}{\sin(kl)}.
    \]
    and the Neumann conditions to obtain a set of $V$ homogeneous equations
    for the variables $f_{i}$ ($i=1,2,\dots V$).\\[0.1cm]
  \item Denote by $\vec{f}$ the vector whose entries are all the $f_{i}$
    variables. Assume that $\sin(kl)\neq 0$ and manipulate the
    linear set of equations you got in the previous section to have the
    following form
    \[
    A\,\vec{f}=\cos(kl)\vec{f}.
    \]
    What is the matrix $A$? Note that this matrix describes the underlying
    discrete graph and this establishes a spectral connection between
    the discrete and the quantum graph.\\[0.1cm]
    % , on which we will elaborate more
    % in the next lectures.\\
  \item Denote by $\left\{ \lambda_{i}\right\} _{i=1}^{V}$ the eigenvalues
    of $A$. Express the $k$-eigenvalues of the quantum graph (remember
    that there are infinitely many of those) in terms of the eigenvalues
    of $A$. Are all the eigenvalues of the quantum graph can be obtained
    in this way? If so, prove it, or otherwise, point on the eigenvalues
    which are not obtained in this way.\\[0.1cm]
  \end{enumerate}

  \textit{Comment}: Further reading on the spectral connection between discrete and quantum
  graphs may be found in \cite{MR813056, Cat_mm97, HarrWey_lmp18, Kuc_wrm04, BelMug_laa13}. The most general derivation of this connection, treating electric and magnetic potentials as well as $\delta$-type vertex conditions appears in \cite{Pan_lmp06}.\\[0.2cm]
\end{itemize}

\section*{acknowledgments}
We thank the referee whose critical remarks and constructive suggestions has lead to a profound improvement
of the manuscript.
RB was supported by ISF (Grant No. 494/14) and Marie Curie Actions (Grant No. PCIG13-GA-2013-618468). \\[0.2cm]

\bibliographystyle{abbrv}
\bibliography{GlobalBib}

\end{document}